\documentstyle[12pt,fleqn,leqno]{article}
\newtheorem{thm}{Theorem}[section]
\newcommand{\be}{\begin{equation}}
\newcommand{\ee}{\end{equation}}
\newcommand{\ba}{\begin{array}}
\newcommand{\ea}{\end{array}}

\renewcommand{\a}{\alpha}
\renewcommand{\b}{\beta}
\renewcommand{\l}{\lambda}
\renewcommand{\t}{\theta}

\renewcommand{\em}{\it}
\newcommand{\bea}{\begin{eqnarray}}
\newcommand{\eea}{\end{eqnarray}}

\newcommand{\Summ}{\sum_{n=0}^\infty}

\begin{document}
\newtheorem{pro}[thm]{Proposition}
\newtheorem{lem}[thm]{Lemma}
\newtheorem{cor}[thm]{Corollary}
\title{\bf $Q$-Hermite Polynomials and Classical Orthogonal
 Polynomials  \thanks{
Research partially supported  by NSF grant DMS 9203659} }
\author{Christian Berg and Mourad E. H.  Ismail}
\date{May 19,1994}

\maketitle
\begin{abstract} We use generating functions to express
orthogonality
 relations in the form of  $q$-beta integrals. The integrand of
such
a  $q$-beta
 integral is then used as  a weight function for a new set of
orthogonal or
 biorthogonal
functions. This method is applied to the continuous  $q$-Hermite
polynomials,
the Al-Salam-Carlitz polynomials, and the polynomials of Szeg\H{o}
and
 leads na\-tu\-rally to the Al-Salam-Chihara polynomials then to
the
 Askey-Wilson polynomials,  the big $q$-Jacobi polynomials and the
 biorthogonal  rational functions
of Al-Salam and Verma, and some recent biorthogonal functions
 of Al-Salam and Ismail.
\end{abstract}

\bigskip
{\bf Running title}: Classical Orthogonal Polynomials.


{\
em Key words and phrases}. Askey-Wilson polynomials,
$q$-orthogonal polynomials, orthogonality relations, $q$-beta
integrals,
$q$-Hermite polynomials.

\bigskip

\setcounter{section}{1}

{\bf 1. Introduction}. The $q$-Hermite polynomials seem
to be at the bottom of a hierarchy
of the classical $q$-orthogonal polynomials, \cite {An:As}.  They
contain no parameters, other than $q$, and one can get them as
special
or limiting cases of other orthogonal polynomials.

The purpose of this work is to show how one can systematically
build the
 classical $q$-orthogonal polynomials from the $q$-Hermite
polynomials
using a simple procedure of attaching generating functions to
measures.

Let $\{p_n(x)\}$ be orthogonal polynomials with respect
to a positive measure $\mu$ with moments of any order and infinite
support 
such that
\be
\int^{\infty}_{-\infty} p_n(x) p_m(x) \, d\mu(x) = \zeta_n
\delta_{m,n}.
\ee
Assume that we know a generating function for $\{p_n(x)\}$, that is
we have
\be
\sum_{n = 0}^\infty p_n(x) t^n/c_n = G(x,t),
\ee
for a suitable numerical sequence  of nonzero elements $\{c_n\}$.
This
 implies that the orthogonality relation (1.1) is  equivalent to
\be \int_{-\infty}^{\infty} G(x, t_1) G(x, t_2) d\mu(x) =
\sum_0^\infty
\zeta_n \frac{(t_1t_2)^n}{c_n^2},
\ee
provided that we can justify the interchange of integration and
sums.

Our idea is to use
\[G(x,t_1) G(x, t_2) \,d\mu(x)\]
 as a new measure, the total mass of which is
given by (1.3), and then look for a system
of functions  (preferably  polynomials) orthogonal or biorthogonal
with
 respect to it. If such a system
is found  one can then repeat the process. It it clear that we
cannot
 indefinitely continue  this process. Things will become too
complicated
at a certain level, and the process will then terminate.

If $\mu$ has compact support it will often be the case that (1.2)
converges
uniformly for $x$ in the support and $|t|$ sufficiently small. In
this
case the justification is obvious.

We mention the following general result with no assumptions about
the support
of $\mu$. For $0<\rho\leq\infty$ we denote by $D(0,\rho)$ the set
of 
$z\in\mbox{\bf C}$ with $|z|<\rho$.

\begin{pro} Assume that (1.1) holds and that the power series
\be
\sum_{n=0}^\infty\frac{\sqrt{\zeta_n}}{c_n}z^n
\ee
has a radius of convergence $\rho$ with $0<\rho\leq\infty$.

(i) Then there is a $\mu$-null set $N\subseteq\mbox{\bf R}$ such
that (1.2)
 converges
 absolutely for $|t|<\rho, x\in \mbox{\bf R}\setminus N$.
Furthermore (1.2)
converges in $L^2(\mu)$ for $|t|<\rho$, and (1.3) holds for
 $|t_1|,|t_2|<\rho$.

(ii) If $\mu$ is indeterminate then (1.2) converges absolutely and
uniformly
 on compact subsets of $\Omega=\mbox{\bf C}\times D(0,\rho)$,
and $G$ is holomorphic in $\Omega$.
\end{pro}

\noindent{\bf Proof}. For $0<r_0<r<\rho$ there exists $C>0$ such
that
$(\sqrt{\zeta_n}/|c_n|)r^n\leq C$ for $n\geq 0$, and we find
\[
\|\sum_{n=0}^N|p_n(x)|\frac{r_0^n}{|c_n|}\|_{L^2(\mu)}\leq
\sum_{n=0}^N\frac{\sqrt{\zeta_n}}{|c_n|}r^n(\frac{r_0}r)^n\leq
C\Summ (\frac{r_0}r)^n <\infty,
\]
which by the monotone convergence theorem implies that
\[
\Summ |p_n(x)|\frac{r_0^n}{|c_n|}\in L^2(\mu),
\]
and in particular the sum is finite for $\mu$-almost all $x$. This
implies
that there is a $\mu$-null set $N\subseteq\mbox{\bf R}$  such that
$\sum p_n(x)(t^n/c_n)$ is absolutely convergent for $|t|<\rho$ and
$x\in
\mbox{\bf R}\setminus N$.

The series (1.2) can be considered as a power series with values in
$L^2(\mu)$,
and by assumption its radius of convergence is $\rho$. It follows
that (1.2)
converges to $G(x,t)$ in $L^2(\mu)$ for $|t|<\rho$, and (1.3) is a
consequence
of Parseval's formula.

 If $\mu$ is indeterminate it is well known that
$\sum|p_n(x)|^2/\zeta_n$
converges uniformly on compact subsets of $\mbox{\bf C}$, cf. \cite
{Ak},
 \cite {Sh:Ta}, and
the assertion follows. \quad $\Box$

\medskip

In order to describe details of our work we will need to introduce
 some notations.
There are  three systems of $q$-Hermite polynomials. Two of them
are
 orthogonal on
  compact subsets of the real line and the third is orthogonal on
the unit circle. The two $q$-Hermite polynomials on the real line
are
 the discrete $q$-Hermite polynomials
$\{H_n(x:q)\}$ and the continuous $q$-Hermite polynomials
$\{H_n(x|q)\}$
 of L. J. Rogers.
 They are generated by
\be
2x H_n(x|q) =  H_{n+1}(x|q) + (1 - q^n)  H_{n-1}(x|q),
\ee
\be
x H_n(x:q) =  H_{n+1}(x:q) + q^{n-1}(1 - q^n)  H_{n-1}(x:q),
\ee
and the initial conditions
\be
H_0(x|q) = H_0(x:q) = 1, \quad H_1(x|q) =  2x,\; H_1(x:q) = x.
\ee
We will describe the $q$-Hermite polynomials on the unit circle
later in
the Introduction. The discrete and continuous $q$-Hermite
polynomials
 have generating functions
\be
\sum_{0}^\infty \frac{H_n(x:q)}{(q;q)_n}t^n =
\frac{(t, -t; q)_\infty}{(xt; q)_\infty},
\ee
and
\be
\sum_{0}^\infty \frac{H_n(x|q)}{(q;q)_n}t^n =
\frac{1}{(te^{i\t}, te^{-i\t}; q)_\infty}, \quad x = \cos \t,
\ee
respectively, where we used the notation in \cite {Ga:Ra} for the
$q$-shifted
 factorials
\be
(a; q)_0 := 1, \quad (a; q)_n := \prod_{k = 1}^n (1 - a
q^{k-1}),\quad n =
1, 2, \cdots, \; or\;  \infty,
\ee
and the multiple $q$-shifted factorials
\be
(a_1, a_2, \cdots, a_k; q)_n := \prod_{j = 1}^k (a_j; q)_n.
\ee
A basic hypergeometric series is
\bea
{}_{r}\phi_s\left(\left. \ba{c} a_1,\ldots,a_{r} \\ b_1,\ldots,b_s
\ea
 \right|\,q,\;z \right) &=& {}_{r}\phi_s(a_1, \ldots, a_{r}; b_1,
\ldots, b_s; q, z) \\
&=& \sum_{n=0}^\infty
\frac{(a_1,\ldots,a_{r};q)_n}{(q,b_1,\ldots,b_s;q)_n}\,z^n
 ((-1)^nq^{n(n-1)/2})^{s+1-r}. \nonumber
\eea

In Section 2 we apply the procedure outlined at the beginning of
the
Introduction to the conti\-nu\-ous $q$-Hermite polynomials for $|q|
< 1$
 and we reach
the Al-Salam-Chihara polynomials in the first step and the second
step takes us
to  the
 Askey-Wilson polynomials. It is worth mentioning
that the Askey-Wilson polynomials are the general classical
orthogonal
polynomials, \cite {An:As}. As a byproduct we get a simple
evaluation
 of the Askey-Wilson $q$-beta integral, \cite {As:Wi}. This seems
to
be the end of the line in this direction. The case $q > 1$ will be
studied in Section 5, see comments below. In Section 3 we apply the
same procedure to the polynomials $\{U_n^{(a)}(x; q)\}$ and
$\{V_n^{(a)}(x; q)\}$ of Al-Salam and Carlitz \cite {Al:Ca}. They
are generated by the recurrences
\be
U_{n+1}^{(a)}(x; q) = [x - (1+a)q^n] U_n^{(a)}(x; q) +aq^{n-1}(1 -
q^n)
 U_{n-1}^{(a)}(x; q), \; n > 0,
\ee
\be
V_{n+1}^{(a)}(x; q) = [x - (1+a)q^{-n}] V_n^{(a)}(x; q)
-aq^{1-2n}(1 - q^n)
 V_{n-1}^{(a)}(x; q), \; n > 0,
\ee
and the initial conditions
\be
 U_0^{(a)}(x; q) =  V_0^{(a)}(x; q) = 1, \; U_1^{(a)}(x; q) =
V_1^{(a)}(x; q)
= x - 1 - a,
\ee
\cite {Al:Ca}, \cite {Ch}.
It is clear that $U_n^{(a)}(x; 1/q) = V_n^{(a)}(x; q)$, so there is
no loss
 of generality in assuming $0 < q < 1$. The $U_n$'s provide a one
parameter
 extension of the discrete $q$-Hermite polynomials when $0 < q < 1$
corresponding to $a=-1$.
In Section 3 we show that
our attachment procedure generates the big $q$-Jacobi polynomials
from the
 $U_n$'s. The big $q$-Jacobi polynomials  were introduced by
Andrews
 and Askey in 1976. The
application of our procedure to the $V_n$'s does not lead
to orthogonal polynomials but to a system of biorthogonal rational
functions of Al-Salam and Verma \cite {Al:Ve}.

The $q$-analogue of Hermite polynomials on the unit circle are the
polynomials
\be
{\cal H}_n(z; q) = \sum_{k=0}^n \frac{(q; q)_n\,(q^{-1/2}z)^k}
{(q; q)_k\,(q; q)_{n-k}}.
\ee
Szeg\H{o} introduced these polynomials in \cite {Sz} to illustrate
his
theory of polynomials orthogonal on the unit circle. Szeg\H{o}
used the Jacobi triple product identity to prove the orthogonality
relation
\be
\frac{1}{2\pi}\int_0^{2\pi}  {\cal H}_m(e^{i\t}; q) \overline{{\cal
H}_n(e^{i\t}; q)} \, (q^{1/2}e^{i\t}, q^{1/2}e^{-i\t}; q)_\infty\,
d\t
 = \frac{(q; q)_n\,q^{-n}}{(q; q)_\infty } \delta_{m,n}.
\ee
In Section 4 we show how generating functions transform (1.16) to
a
$q$-beta integral of Ramanujan. This  explains the origin of the
 biorthogonal polynomials of Pastro \cite {Pa} and the ${}_4\phi_3$
biorthogonal rational functions of Al-Salam and Ismail \cite
{Al:Is}.
                                                
In section 5 we consider the $q$-Hermite polynomials for $q>1$.
They are
orthogonal on the imaginary axis. For $0<q<1$ we put
$h_n(x|q)=(-i)^nH_n(ix|1/q)$,
and $\{h_n(x|q)\}$ are called the $q^{-1}$-Hermite polynomials.
They correspond
to an indeterminate moment problem considered in detail in \cite
{Is:Ma2}.
Using a $q$-analogue of the Mehler formula for these polynomials we
derive
an analogue of the Askey-Wilson integral valid for all the
solutions to the
indeterminate moment problem. Our derivation, which is different
from the one
in \cite {Is:Ma2}, is based on Parseval's formula.

The attachment procedure for the $q^{-1}$-Hermite polynomials leads
to a
special case of the Al-Salam-Chihara polynomials corresponding to
$q>1$, more
precisely to the polynomials
\be
u_n(x; t_1, t_2) = v_n(-2x; q, -(t_1+t_2)/q, t_1t_2q^{-2}, -1),
\ee
cf. \cite {As:Is2}. We prove that for any positive orthogonality
measure
$\mu$ for the $q^{-1}$-Hermite polynomials
\be
d\nu_{\mu}(\sinh\xi; t_1, t_2) := \frac{(-t_1e^{\xi},t_1e^{-\xi},
-t_2e^{\xi},t_2e^{-\xi}; q)_\infty}
{(-t_1t_2/q; q)_\infty} d\mu(\sinh\xi),
\ee
is an orthogonality measure for $\{u_n\}$.

The attachment procedure applied to $\{u_n\}$ leads to the
biorthogonal
rational functions    
\be
\varphi_n (\sinh\,\xi ;t_1,t_2,t_3,t_4) :=\,_4\phi_3
\left(\begin{array}{cc} q^{-n},-t_1t_2q^{n-2},-
t_1t_3/q,-t_1t_4/q\\
 -t_1e^\xi ,t_1e^{-\xi}  ,t_1t_2t_3t_4q^{-3}\end{array}\;\vrule\;\;
  q,q\right).
\ee
of Ismail and Masson \cite {Is:Ma2} in the special case
$t_3=t_4=0$.
\bigskip

\setcounter{equation}{0}

\setcounter{section}{2}

{\bf  2. The Continuous $q$-Hermite Ladder}. Here we assume
$-1<q<1$.
The orthogonality relation
for the continuous $q$-Hermite polynomials is
\be
\int_0^\pi H_m(\cos \t|q)H_n(\cos \t|q) (e^{2i\t}, e^{-2i\t};
q)_\infty d\t =
\frac{2\pi (q; q)_n}{(q; q)_\infty} \delta_{m,n}.
\ee
The series in (1.9) converges for $|t|<1$ uniformly in $\t\in
[0,\pi]$.
Thus (1.2), (1.3) and the generating function (1.9) imply
\be
\int_0^\pi \frac{(e^{2i\t}, e^{-2i\t}; q)_\infty}
{(t_1e^{i\t}, t_1e^{-i\t}, t_2e^{i\t}, t_2e^{-i\t}; q)_\infty} d\t
= \frac {2\pi}{(q, t_1t_2; q)_\infty}, \quad |t_1|, |t_2|<1,
\ee
where we used the $q$-binomial theorem \cite [(II.3)]{Ga:Ra}
\be
\sum_0^\infty \frac{(a; q)_n}{(q; q)_n} z^n =
\frac{(az; q)_\infty}{(z; q)_\infty},
\ee
with $a = 0$.

The next step is to find polynomials $\{p_n(x)\}$ orthogonal with
respect to
the weight function
\be
w_1(x; t_1, t_2) := \frac{(e^{2i\t}, e^{-2i\t}; q)_\infty}
{(t_1e^{i\t}, t_1e^{-i\t}, t_2e^{i\t}, t_2e^{-i\t}; q)_\infty}
  \frac{1}{\sqrt{1-x^2}}, \quad x = \cos \t,
\ee
which is positive for $t_1,t_2\in \; (-1,1).$
Here we follow a clever technique of attachment which was used by
Askey and
 Andrews and by Askey and Wilson in \cite{As:Wi}.
Write $\{p_n(x)\}$ in the form
\be
p_n(x) = \sum_{k=0}^n \frac{(q^{-n}, t_1e^{i\t}, t_1e^{-i\t};
q)_k}{(q; q)_k}a_{n,k},
\ee
then determine $a_{n,k}$ such that $p_n(x)$ is orthogonal to
$(t_2e^{i\t},
t_2e^{-i\t}; q)_j$, $j = 0$, $1$, $\cdots$, $n-1$. Note that
$(ae^{i\t}, ae^{-i\t}; q)_k$ is a polynomial in $x$ of degree $k$,
since
\bea
(ae^{i\t}, ae^{-i\t}; q)_k = \prod_{j = 0}^{k-1}(1 - 2axq^j +
a^2q^{2j}).
\eea
The reason for choosing the bases $\{(t_1e^{i\t},
t_1e^{-i\t}; q)_k\}$ and $\{(t_2e^{i\t},
t_2e^{-i\t}; q)_j\}$ is that they attach nicely to the weight
function
and (2.2)  enables us  to integrate  $(t_1e^{i\t},
t_1e^{-i\t}; q)_k (t_2e^{i\t},
t_2e^{-i\t}; q)_j$ against  the weight function $w_1(x; t_1, t_2)$.
Indeed
\bea
(t_1e^{i\t},
t_1e^{-i\t}; q)_k (t_2e^{i\t},
t_2e^{-i\t}; q)_j w_1(x; t_1, t_2) = w_1(x; t_1q^k, t_2q^j).
\nonumber
\eea
Therefore
\bea
\int_{-1}^1 (t_2e^{i\t}, t_2e^{-i\t}; q)_j p_n(x) w_1(x; t_1, t_2)
dx
\nonumber
\eea
\bea
= \sum_{k = 0}^n \frac{(q^{-n}; q)_k}{(q; q)_k} a_{n,k}\int_0^{\pi}
\frac{(e^{2i\t}, e^{-2i\t}; q)_\infty \; d\t}{(t_1q^ke^{i\t},
t_1q^ke^{-i\t},
t_2q^je^{i\t}, t_2q^je^{-i\t}; q)_\infty}
\nonumber
\eea
\bea
= \frac{2\pi}{(q; q)_\infty} \sum_{k = 0}^n \frac{(q^{-n}; q)_k
a_{n,k}}
{(q; q)_k (t_1t_2q^{k+j}; q)_\infty}
\nonumber
\eea
\bea
= \frac{2\pi}{(q, t_1t_2q^j; q)_\infty} \sum_{k = 0}^n
\frac{(q^{-n}, t_1t_2q^j; q)_k}{(q; q)_k} a_{n,k}.
\nonumber
\eea
At this stage we look for $a_{n,k}$ as a quotient of products of
 $q$-shifted factorials in order to make the above sum vanish for
$0 \le j <n$. The $q$-Chu-Vandermonde sum \cite [(II.6)]{Ga:Ra}
\be
{}_2\phi_1(q^{-n}, a; c; q, q) = \frac{(c/a; q)_n}{(c;q)_n}a^n
\ee
suggests
\bea
a_{n,k} = q^k/(t_1t_2;q)_k. \nonumber
\eea
Therefore
\bea
\int_{-1}^1 (t_2e^{i\t}, t_2e^{-i\t}; q)_j p_n(x) w_1(x; t_1, t_2)
dx
 = \frac{2\pi(q^{-j}; q)_n}{(q, t_1t_2q^j; q)_\infty\,
(t_1t_2;q)_n}\,(t_1t_2q^j)^n. \nonumber
\eea
It follows from (2.5) and (2.6)  that the coefficient of $x^n$ in
$p_n(x)$
is
\be
 (-2t_1)^nq^{n(n+1)/2}(q^{-n}; q)_n/(q, t_1t_2;
q)_n=(2t_1)^n/(t_1t_2;q)_n.
\nonumber
\ee
This lead to the orthogonality relation
\be
\int_{-1}^1 p_m(x) p_n(x) w_1(x; t_1, t_2) dx
= \frac{2\pi (q; q)_n t_1^{2n}}{(q, t_1t_2; q)_\infty (t_1t_2;
q)_n}
\delta_{m,n}.
\ee
Furthermore the polynomials are given by
\be
p_n(x) = {}_{3}\phi_2\left(\left. \ba{c} q^{-n}, t_1e^{i\t},
t_1e^{-i\t} \\
t_1t_2, 0 \ea \right|\,q,\;q \right).
\ee
The polynomials we have just found are the Al-Salam-Chihara
polynomials
and were first identified by  W. Al-Salam and  T. Chihara  \cite
{Al:Ch}.
Their weight function was given  in \cite {As:Is2} and \cite
{As:Wi}.

Observe that the orthogonality relation (2.9) and the uniqueness of
the
 polynomials orthogonal with respect to a positive measure show
that
 $t_1^{-n} p_n(x)$ is symmetric in $t_1$ and $t_2$. This gives the
known
 transformation
\be
{}_{3}\phi_2\left(\left. \ba{c} q^{-n}, t_1e^{i\t}, t_1e^{-i\t} \\
t_1t_2, 0 \ea \right|\,q,\;q \right) =
(t_1/t_2)^n{}_{3}\phi_2\left(\left. \ba{c} q^{-n}, t_2e^{i\t},
t_2e^{-i\t} \\
t_1t_2, 0 \ea \right|\,q,\;q \right)
\ee
as a byproduct of our analysis.

Our next task is to repeat the process with the Al-Salam-Chihara
polynomials
 as our starting point.  The representation (2.10) needs to be
transformed
 to a form
more amenable to generating functions. This can be done using an
idea of
 Ismail and Wilson \cite {Is:Wi}. First write the $_3\phi_2$ as a
sum over $k$
 then replace $k$ by $n-k$. Applying
\be
(a; q)_{n-k} = \frac{(a; q)_n}{(q^{1-n}/a;q)_k}\,
(-q/a)^kq^{-kn+k(k-1)/2}
\ee
we obtain
\bea
p_n(x) = \frac{(t_1e^{i\t}, t_1e^{-i\t};q)_n}{(t_1t_2;
q)_n}q^{-n(n-1)/2}(-1)^n
\sum_{k=0}^n \frac{(-t_2/t_1)^k\, (q^{-n}, q^{1-n}/t_1t_2; q)_k}
{(q,  q^{1-n}e^{i\t}/t_1, q^{1-n}e^{-i\t}/t_1; q)_k}q^{k(k+1)/2}.
\eea
Applying the $q$-analogue of the Pfaff-Kummer transformation
 \cite [(III.4)]{Ga:Ra}
\be
\Summ \frac{(A, C/B; q)_n}{(q, C, Az;q)_n}q^{n(n-1)/2}(-Bz)^n
= \frac{(z; q)_\infty}{(Az; q)_\infty} {}_2\phi_1(A, B; C; q, z),
\ee
with
\bea
 A = q^{-n},\; B = t_2 e^{i\t}, \; C = q^{1-n} e^{i\t}/t_1, \; z =
qe^{-i\t}/t_1
\nonumber
\eea
to (2.13) we obtain the representation
\bea
p_n(x) = \frac{(t_1e^{-i\t}; q)_n t_1^ne^{in\t}}{(t_1t_2; q)_n}
\, {}_{2}\phi_1\left(\left. \ba{c} q^{-n}, t_2e^{i\t} \\
q^{1-n}e^{i\t}/t_1 \ea \right|\,q,\;qe^{-i\t}/t_1 \right).  
\nonumber
\eea
Using (2.12) we express a multiple of $p_n$ as a Cauchy product of
two
 sequences. The result is
\bea
p_n(x) = \frac{(q; q)_n t_1^n}{(t_1t_2; q)_n} \sum_{k = 0}^n
\frac{(t_2e^{i\t}; q)_k}{(q; q)_k} e^{-ik\t}\;
\frac{(t_1e^{-i\t}; q)_{n-k}}{(q;q)_{n-k}} e^{i(n-k)\t}. \nonumber
\eea
This and the $q$-binomial theorem (2.3) establish the generating
function
\be
\Summ \frac{(t_1t_2; q)_n}{(q; q)_n}\, p_n(x) \,(t/t_1)^n =
\frac{(tt_1, tt_2;q)_\infty}{(te^{-i\t}, te^{i\t}; q)_\infty}
\ee
The orthogonality relation (2.9) and the generating function (2.15)
 imply the Askey-Wilson $q$-beta integral, \cite {As:Wi}, \cite
{Ga:Ra}
\be
\int_{0}^\pi\frac{(e^{2i\t}, e^{-2i\t}; q)_\infty}
{\prod_{j=1}^4 (t_je^{i\t}, t_j e^{-i\t}; q)_\infty} d\t
= \frac{2\pi\, (t_1t_2t_3t_4; q)_\infty}{(q; q)_\infty
\prod_{1\le j <k \le 4} (t_jt_k; q)_\infty}.
\ee
The polynomials orthogonal with respect to the weight function
 whose total mass is given by (2.16) are the Askey-Wilson
polynomials.
Their explicit representation and orthogonality follow from (2.16)
and the
$q$-analogue of the Pfaff-Saalsch\"{u}tz theorem, \cite
[(II.12)]{Ga:Ra}.
The details of this  calculation are in \cite {As:Wi}. The
polynomials are
 \be
p_n(x; t_1, t_2, t_3, t_4|q)  = t_1^{-n}(t_1t_2, t_1t_3,
t_1t_4;q)_n\;
{}_{4}\phi_3\left(\left. \ba{c} q^{-n},
t_1t_2t_3t_4q^{n-1}, t_1e^{i\t}, t_1e^{-i\t} \\
t_1t_2,\; t_1t_3, \;t_1t_4 \ea \right|\,q,\;q \right).
\ee

The orthogonality relation of the Askey-Wilson polynomials is \cite
[(2.3)-(2.5)]{As:Wi}
\be
\int_0^\pi p_m(\cos \t; t_1,t_2, t_3, t_4|q)
p_n(\cos \t; t_1,t_2, t_3, t_4|q) w(\cos \t;  t_1,t_2, t_3, t_4)
d\t
\ee
\bea
\qquad = \frac{2 \pi\, (t_1t_2t_3t_4q^{2n};q)_\infty\,
(t_1t_2 t_3 t_4q^{n-1};q )_n}{(q^{n+1}; q)_\infty\, \prod_{1 \le j
< k \le 4}
(t_jt_kq^n;q)_\infty}\delta_{m,n}, \nonumber
\eea
for $\max\{|t_1|, |t_2, |t_3|, |t_4|\} < 1$ and the weight function
is given by
\be
 w(\cos \t;  t_1,t_2, t_3, t_4) = \frac{(e^{2i\t}, e^{-2i\t};
q)_\infty}
{\prod_{j=1}^4 (t_je^{i\t}, t_j e^{-i\t}; q)_\infty}.
\ee
Observe that the weight function in (2.19) and the right-hand side
of
(2.18) are symmetric functions of $t_1, t_2, t_3, t_4$. The weight
function
in (2.18) is positive when $\max\{|t_1|, |t_2, |t_3|, |t_4|\} < 1$
and the
 uniqueness of the polynomials orthogonal with respect to a
positive
 measure shows that the Askey-Wilson polynomials are symmetric in
the
four parameters $t_1, t_2, t_3, t_4$. This symmetry is the Sears
transformation
 \cite [(III.15)]{Ga:Ra}, a fundamental transformation
in the theory of basic hypergeometric functions. The Sears
transformation
may be stated in the form
\be
{}_{4}\phi_3\left(\left. \ba{c} q^{-n}, a, b, c \\
d,\; e, \;f \ea \right|\,q,\;q \right) =
\left(\frac{bc}{d}\right)^n
\frac{(\frac{de}{bc}, \frac{df}{bc})_n}{(e, f)_n}\;
{}_{4}\phi_3\left(\left. \ba{c} q^{-n}, a, d/b, d/c \\
d,\; \frac{de}{bc}, \frac{df}{bc} \ea \right|\,q,\;q \right),
\ee
where $abc = def q^{n-1} $.

Ismail and Wilson  \cite {Is:Wi} used the Sears transformation to
establish
 the generating function
\be
\Summ \frac{p_n(\cos\t; t_1, t_2, t_3, t_4|q)}{(q, t_1t_2, t_3t_4;
q)_n}\,
t^n = {}_{2}\phi_1\left(\left. \ba{c} t_1e^{i\t}, t_2e^{i\t} \\
t_1t_2 \ea \right|\,q,\;te^{-i\t} \right)\,
{}_{2}\phi_1\left(\left. \ba{c} t_3e^{-i\t}, t_4e^{-i\t} \\
t_3t_4 \ea \right|\,q,\;te^{i\t} \right).
\ee
Thus (2.18) leads to the evaluation of the following integral
\be
\int_0^\pi \prod_{j=5}^6 {}_{2}\phi_1\left(\left. \ba{c}
t_1e^{i\t}, t_2e^{i\t} \\
t_1t_2 \ea \right|\,q,\;t_je^{-i\t} \right)\,
{}_{2}\phi_1\left(\left. \ba{c} t_3e^{-i\t}, t_4e^{-i\t} \\
t_3t_4 \ea \right|\,q,\;t_je^{i\t} \right)
\ee
\bea
\quad \quad  \times \frac{(e^{2i\t}, e^{-2i\t}; q)_\infty}
{\prod_{j=1}^4 (t_je^{i\t}, t_j e^{-i\t}; q)_\infty} \, d\t
\nonumber
\eea
\bea
= \frac{2\pi\, (t_1t_2t_3t_4; q)_\infty}{(q; q)_\infty\,
 \prod_{1 \le j < k \le 4} (t_jt_k; q)_\infty}\,
{}_{6}\phi_5\left(\left. \ba{c} \sqrt{t_1t_2t_3t_4/q}, -
\sqrt{t_1t_2t_3t_4/q},
t_1t_3, t_1t_4, t_2t_3, t_2t_4 \\
\sqrt{t_1t_2t_3t_4q}, - \sqrt{t_1t_2t_3t_4q}, t_1t_2, t_3t_4,
t_1t_2t_3t_4/q \ea \right|\,q,\;t_5t_6 \right), \nonumber
\eea
valid for $\max\{|t_1|, |t_2|, |t_3|,|t_4|,|t_5|,|t_6|\} < 1$.

\bigskip

\setcounter{equation}{0}

\setcounter{section}{3}

{\bf  3. The Discrete $q$-Hermite Ladder}. Here we assume $0<q<1$.
Instead of using the
discrete $q$-Hermite polynomials directly we will use the
Al-Salam-Carlitz $q$-polynomials which are a one parameter
generalization of the discrete $q$-Hermite polynomials. The
 Al-Salam-Carlitz polynomials $\{U_n^{(a)}(x; q)\}$  have the
generating function \cite {Al:Ca}, \cite {Ch}
\be
G(x; t) := \Summ U_n^{(a)}(x; q)\frac{t^n}{(q; q)_n}
= \frac{(t, at;q)_\infty}{(tx; q)_\infty}, \quad a < 0,\; 0 < q <
1,
\ee
and  satisfy the the orthogonality relation
\be
\int_{-\infty}^\infty U_m^{(a)}(x; q)U_n^{(a)}(x; q)\,
d\mu^{(a)}(x)
= (-a)^nq^{n(n-1)/2}(q; q)_n\delta_{m,n},
\ee
with $\mu^{(a)}$ a discrete probability measure on $[a, 1]$  given
by
\be
\mu^{(a)} =  \Summ\left[\frac{q^n}{(q, q/a; q)_n(a; q)_\infty}
 \varepsilon_{q^n} + \frac{q^n}{(q, aq; q)_n(1/a; q)_\infty}
 \varepsilon_{aq^n}\right].
\ee
In (3.3) $\varepsilon_y$ denotes a unit mass supported at $y$.
 The form of the orthogonality relation (3.2)-(3.3) given in \cite
{Al:Ca}
and \cite {Ch} contained a complicated form of a normalization
constant.
The value of the constant was simplified in \cite {Is}.
Since the radius of convergence of (1.4) is $\rho=\infty$ we can
apply
Proposition 1.1 and get
\bea
 \int_{-\infty}^\infty G(x; t_1)\, G(x; t_2)\,  d\mu^{(a)}(x) =
\Summ
\frac{(-at_1t_2)^n}{(q; q)_n} q^{n(n-1)/2} = (at_1t_2;
q)_\infty,\quad
t_1,t_2\in{\bf C},
\eea
where the last expression follows by Euler's theorem \cite
[(II.2)]{Ga:Ra}
\be
\Summ z^nq^{n(n-1)/2}/(q; q)_n = (-z; q)_\infty.
\ee
This establishes the integral
\be
\int_{-\infty}^\infty \frac{d\mu^{(a)}(x)}{(xt_1, xt_2; q)_\infty}
= \frac{(at_1t_2; q)_\infty}{(t_1, t_2, at_1, at_2; q)_\infty}.
\ee

When we substitute for $\mu^{(a)}$ from (3.3) in (3.4) or (3.6) we
discover the nonterminating Chu-Vandermonde sum,
\cite[(II.23)]{Ga:Ra},
\bea
\frac{(Aq/C, Bq/C;q)_\infty}{(q/C;
q)_\infty}{}_{2}\phi_1\left(\left. \ba{c}
A, B \\
C \ea \right|\,q,\;q \right) + \frac{(A, B;q)_\infty}{(C/q;
q)_\infty}
{}_{2}\phi_1\left(\left. \ba{c}
Aq/C, Bq/C \\
q^2/C \ea \right|\,q,\;q \right)  \\
= (ABq/C;q)_\infty \nonumber
\eea

We now restrict the attention to $t_1,t_2\in \,(a^{-1},1)$ in the
case of which $1/(xt_1,xt_2; q)_\infty$ is a positive weight
function on
$[a,1]$.
The next step is to find polynomials orthogonal with respect to
$d\mu^{(a)}(x)/(xt_1, xt_2; q)_\infty$. Define $P_n(x)$ by
\be
P_n(x) = \sum_{k=0}^n \frac{(q^{-n}, xt_1; q)_k}{(q; q)_k}\, q^k
a_{n,k}
\ee
where $a_{n,k}$ will be chosen later. Using (3.6) it is easy to see
that
\bea
\int_{-\infty}^\infty P_n(x)
\frac{(xt_2; q)_m}{(xt_1, xt_2; q)_\infty} \,d\mu^{(a)}(x)=
\sum_{k=0}^n
\frac{(q^{-n}; q)_k}{(q; q)_k}\, q^k a_{n,k}
\frac{(at_1t_2q^{k+m}; q)_\infty}{(t_1q^k,  at_1q^k, t_2q^m,
at_2q^m; q)_\infty}
 \nonumber
\eea
\bea
= \frac{(at_1t_2q^{m}; q)_\infty}{(t_1,  at_1, t_2q^m, at_2q^m;
q)_\infty}
 \sum_{k=0}^n \frac{(q^{-n}, t_1, at_1;q)_k}{(q, at_1t_2q^m; q)_k}
a_{n,k}q^k.
\nonumber
\eea
The choice $a_{n,k}= (\l;q)_k/(t_1, at_1; q)_k$ allows us to apply
 the $q$-Chu-Vandermonde sum (2.7). The choice $\l =
at_1t_2q^{n-1}$  leads to
\be
\int_{-\infty}^\infty P_n(x)
\frac{(xt_2; q)_m}{(xt_1, xt_2; q)_\infty}\, d\mu^{(a)}(x)
= \frac{(at_1t_2q^m; q)_\infty\, (q^{m+1-n}; q)_n\,
(at_1t_2q^{n-1})^n}
{(t_1, at_1, t_2q^m, at_2q^m; q)_\infty\, (at_1t_2q^m;q)_n}.
\ee
The right-hand side of (3.9) vanishes for $0 \le m < n$.
The coefficient of $x^n$ in $P_n(x)$ is
\[
\frac{(q^{-n}, at_1t_2q^{n-1}; q)_n}{(q, t_1, at_1;
q)_n}(-t_1)^nq^{n(n+1)/2}
=\frac{(at_1t_2q^{n-1}; q)_n}{(t_1,at_1; q)_n}t_1^n.
\]
Therefore
\be
P_n(x) = \varphi_n(x;  a, t_1, t_2) =  {}_{3}\phi_2\left(\left.
\ba{c} q^{-n}, at_1t_2q^{n-1}, xt_1 \\
t_1, at_1 \ea \right|\,q,\;q \right),
\ee
satisfies the orthogonality relation
\be
\int_{-\infty}^\infty \varphi_m(x; a, t_1, t_2) \varphi_n(x; a,
t_1, t_2)
\frac{d\mu^{(a)}(x)}{(xt_1, xt_2; q)_\infty}
    \ee
\bea \qquad= \frac{(q, t_2, at_2, at_1t_2q^{n-1}; q)_n \,
 (at_1t_2q^{2n}; q)_\infty}
{(t_1, at_1, t_2, at_2; q)_\infty\, (t_1, at_1; q)_n}\;
(-at_1^2)^n\,
q^{n(n-1)/2} \,
\delta_{m,n}. \nonumber
\eea
The polynomials $\{\varphi_n(x;  a, t_1, t_2)\}$ are the big
$q$-Jacobi polynomials of Andrews and Askey \cite {An:As} in
 a different normalization. The Andrews-Askey normalization is
\be
P_n(x; \a, \b, \gamma:q) =  {}_{3}\phi_2\left(\left. \ba{c} q^{-n},
 \a \b q^{n+1}, x \\ \a q, \gamma q \ea \right|\,q,\;q \right).
\ee

Note that we may rewrite the  orthogonality relation (3.11) in the
form
\bea
\int_{-\infty}^\infty t_1^{-m}(t_1, at_1; q)_m\varphi_m(x; a, t_1,
t_2)
\;  t_1^{-n} (t_1, at_1; q)_n
\varphi_n(x; a, t_1, t_2)
\frac{d\mu^{(a)}(x)}{(xt_1, xt_2; q)_\infty}
\eea
\bea \qquad= \frac{(q,t_1, at_1, t_2, at_2, at_1t_2q^{n-1}; q)_n \,
(at_1t_2q^{2n}; q)_\infty}
{(t_1, at_1, t_2, at_2; q)_\infty}\; (-a)^n\,
q^{n(n-1)/2} \,
\delta_{m,n}. \nonumber
\eea
Since $d\mu^{(a)}(x)/(xt_1, xt_2; q)_\infty$ and the right-hand
side of (3.13) are symmetric in $t_1$ and $t_2$ then
\bea
t_1^{-n} (t_1, at_1; q)_n \, \varphi_n(x; a, t_1, t_2) \nonumber
\eea
must be symmetric in $t_1$ and $t_2$. This gives the ${}_3\phi_2$
transformation
\bea
 {}_{3}\phi_2\left(\left. \ba{c} q^{-n}, at_1t_2q^{n-1}, xt_1 \\
t_1, at_1 \ea \right|\,q,\;q \right)
= \frac{t_1^{n} (t_2, at_2; q)_n}{t_2^{n} (t_1, at_1; q)_n}
{}_{3}\phi_2\left(\left. \ba{c} q^{-n}, at_1t_2q^{n-1}, xt_2 \\
t_2, at_2 \ea \right|\,q,\;q \right).
\eea

We now consider the polynomials $\{V_n^{(a)}(x; q)\}$ and restrict
the
parameters to $0<a,\;0<q<1,$ in which case they are orthogonal with
respect
 to a positive measure, cf.  \cite [VI.10]{Ch}. The corresponding
moment
 problem  is determinate if and only if $0<a\leq q$ or $1/q\leq a$.
In the
first case the unique solution is
\be
m^{(a)}=(aq;q)_\infty\sum_{n=0}^\infty
\frac{a^nq^{n^2}}{(q,aq;q)_n}\varepsilon_{q^{-n}},
\ee
and in the second case it is
\be
\sigma^{(a)}=(q/a;q)_\infty\sum_{n=0}^\infty
\frac{a^{-n}q^{n^2}}{(q,q/a;q)_n}\varepsilon_{aq^{-n}},
\ee
cf. \cite{Be:Va}. The total mass of these measures was evaluated to
1 in
\cite {Is}.

If $q<a<1/q$ the problem is indeterminate and both measures are
solutions.
In \cite{Be:Va} the following one-parameter family of solutions
with an
analytic density was found
\be
\nu(x;a,q,\gamma)=\frac{\gamma|a-1|(q,aq,q/a;q)_\infty}
{\pi a[(x/a;q)_{\infty}^2+\gamma^2(x;q)_{\infty}^2]},\quad
\gamma>0.
\ee
In the above
$a=1$ has to be excluded. For a similar formula when $a=1$ see
\cite{Be:Va}.

If $\mu$ is one of the solutions of the moment problem we have the
 orthogonality relation
\be
\int_{-\infty}^\infty V_m^{(a)}(x; q)V_n^{(a)}(x; q) d\mu(x)
 =a^nq^{-n^2}(q; q)_n \delta_{m,n}.
\ee
The polynomials have the
generating function \cite {Al:Ca}, \cite {Ch}
\be
V(x; t) := \Summ V_n^{(a)}(x; q) \frac{q^{n(n-1)/2}}{(q;q)_n}
(-t)^n =
\frac{(xt;q)_\infty}{(t, at; q)_\infty},\quad
|t|<\mbox{min$\;$}(1,1/a).
\ee

The power series (1.4) has the radius of convergence $\sqrt{q/a}$,
and
 therefore (1.3) becomes
\bea
\int_{-\infty}^\infty \frac{(xt_1, xt_2; q)_\infty\; d\mu(x)}
{(t_1, at_1, t_2, at_2; q)_\infty} &= & \int_{-\infty}^\infty V(x,
t_1)\,V(x, t_2) d\mu(x)
 \nonumber \\
&=& \Summ \frac{(at_1t_2/q)^n}{(q; q)_n}  \nonumber \\
&=& \frac{1}{(at_1t_2/q; q)_\infty},\quad |t_1|,|t_2|<\sqrt{q/a} . 
\nonumber
\eea              

This  identity with $\mu=m^{(a)}$ or $\mu=\sigma^{(a)}$ is nothing
but the
 $q$-analogue of the Gauss theorem,
\bea
{}_2\phi_1(a, b; c; q, c/ab)= \frac{(c/a, c/b; q)_\infty}{(c, c/ab;
q)_\infty}
\eea
\cite[(II.8)]{Ga:Ra}.

Specializing  to the density (3.17) we get

\be
\int_{-\infty}^\infty \frac{(xt_1, xt_2; q)_\infty \; dx}
{(x/a;q)^2_\infty+\gamma^2(x;q)^2_\infty} =  \frac {\pi a(t_1,at_1,
t_2,
 at_2; q)_\infty}{|a-1|\gamma (q,aq,q/a, at_1t_2/q; q)_\infty},
\ee
valid for $q<a<1/q, a\neq 1,\gamma >0$.

We now seek polynomials  or rational functions that are orthogonal
with
respect to the measure
\be
d\nu(x)=(xt_1,xt_2;q)_\infty\;d\mu(x),
\ee
where $\mu$ satisfies (3.18).
It is clear that we can integrate
$1/[(xt_1; q)_k(xt_2;q)_j]$ with respect to the measure $\nu$.
Set
\be
\psi_n(x; a, t_1, t_2) := \sum_{k=0}^n \frac{(q^{-n}; q)_k}{(q;
q)_k}
\frac{q^ka_{n,k}}{(xt_1; q)_k}.
\ee

The rest of the analysis is similar to our treatment of the
$U_n$'s.
We get
\[
\int_{-\infty}^\infty \frac{\psi_n(x; a, t_1,
t_2)}{(xt_2;q)_m}\;d\nu(x)=
\sum_{k=0}^n\frac{(q^{-n};q)_k}{(q;q)_k}q^ka_{n,k}\int_{-\infty}^
\infty (xt_1q^k,
xt_2q^m;q)_\infty\;d\mu(x),
\]

and if we choose $a_{n,k} = (t_1, at_1; q)_k/(at_1t_2/q; q)_k$
the above expression is equal to
\bea
\frac{(t_1,at_1,t_2q^m,at_2q^m;q)_\infty}{(at_1t_2q^{m-1};q)_\infty}
\, {}_{2}\phi_1\left(\left. \ba{c} q^{-n}, at_1t_2q^{m-1} \\
at_1t_2/q  \ea \right|\,q,\;q \right),\nonumber
\eea
\bea
\quad=\frac{(t_1,at_1,t_2q^m,at_2q^m;q)_\infty (q^{-m};q)_n }
{(at_1t_2q^{m-1};q)_\infty
(at_1t_2/q;q)_n}(at_1t_2q^{m-1})^n,\nonumber
\eea
which is  0 for $m<n$. We have used the Chu-Vandermonde sum (2.7).
Since $\nu$ is symmetric in $t_1,t_2$ this  leads to the
biorthogonality
 relation
\be
\int_{-\infty}^\infty  \psi_m(x; a, t_2, t_1)\psi_n(x; a, t_1,
t_2)\;d\nu(x)
 = \frac{(t_1, at_1, t_2, at_2; q)_\infty (q;q)_n}
{(at_1t_2/q; q)_\infty (at_1t_2/q; q)_n} (at_1t_2/q)^n\delta_{m,n}.
\ee

The $\psi_n$'s are  given by
\bea
\psi_n(x; a, t_1, t_2) = \,_3\phi_2\left(\begin{array}{cc} q^{-n},
t_1, at_1
\\ xt_1, at_1t_2/q \end{array}\;\;\vrule\;\; q,q\right).
\eea
They are essentially the rational functions studied by Al-Salam and
Verma in \cite {Al:Ve}. Al-Salam and Verma  used the notation
 \be
R_n (x;\alpha ,\beta ,\gamma ,\delta ;q) =
\,_3\phi_2\left(\begin{array}{cc} \beta ,\alpha\gamma /\delta
,q^{-n}\\ \beta\gamma /q,\alpha qx\end{array}\;\;\vrule\;\;
q,q\right).
\ee
The translation between the two notations is
\be
\psi_n(x; a, t_1, t_2) = R_n(\b xq^{-1}/\a; \a, \b , \gamma, \delta
;q),
\ee
with
\be
t_1 = \b, \quad t_2 = \b \delta/q\a, \quad a = \a\gamma/\b\delta.
\ee
Note that $R_n$ has only three free variables since one of the
parameters
 $\a, \b, \gamma, \delta$ can be absorbed by  scaling the
independent variable.

\bigskip

\setcounter{equation}{0}

\setcounter{section}{4}

{\bf  4. The Szeg\H{o} Ladder}. As already mentioned in the
 introduction Szeg\H{o}
\cite {Sz} used the Jacobi triple product identity to prove (1.17).
The explicit form (1.16) and the $q$-binomial theorem (2.3) give
\be
{\cal H}(z, t) := \Summ {\cal H}_n(z; q) \frac{t^n}{(q;q)_n} =
 1/(t, tzq^{-1/2}; q)_\infty.
\ee
Thus (1.17) and (4.1) imply the Ramanujan $q$-beta integral
\be\frac{1}{2\pi i} \int_{|z| = 1} \frac{(q^{1/2}z, q^{1/2}/z;
q)_\infty}
{(t_1q^{-1/2}z, t_2q^{-1/2}/z; q)_\infty} \, \frac{dz}{z} =
\frac{(t_1, t_2; q)_\infty}{(q, t_1t_2/q; q)_\infty},
\ee
for $|t_1| < q^{1/2},\, |t_2| < q ^{1/2}$, since (1.4) has radius
of
 convergence $q^{1/2}$.
Putting
\be
\Omega(z)=\frac{(q,t_1t_2q,q^{1/2}z,q^{1/2}/z;q)_\infty}
{(t_1q,t_2q,t_1q^{1/2}z,t_2q^{1/2}/z;q)_\infty}
\ee
and applying the attachment technique
 to (4.2) we find that the polynomials
\be
\tilde p_n(z,t_1,t_2) :=
{}_{3}\phi_2\left(\left. \ba{c} q^{-n}, t_1q^{1/2}z, t_1q \\
0,\;  t_1t_2q \ea \right|\,q,\;q \right).
\ee
satisfy the biorthogonality relation
\be
\frac{1}{2\pi i}\int_{|z| = 1} \tilde p_m(z, t_1, t_2) \overline
{\tilde p_n(z, \overline{t_2},\overline{t_1})}\,\Omega(z)
\, \frac{dz}{z} = \frac{(q; q)_n}{(t_1t_2q;
q)_n}(t_1t_2q)^n\delta_{m,n}.
\ee
Using the transformation \cite [(III.7)] {Ga:Ra} we see that
\be
\tilde
p_n(z,t_1,t_2)=\frac{(q;q)_n}{(t_1t_2q;q)_n}(t_1q)^np_n(z,t_1,t_2),
\ee
where
\be
p_n(z, a, b) = \frac{(b; q)_n}{(q; q)_n} {}_2\phi_1(q^{-n}, aq;
q^{1-n}/b; q, q^{1/2}z/b) = \sum_{k = 0}^n \frac{(aq;
q)_k\,(b;q)_{n-k}}
{(q; q)_k\,(q; q)_{n-k}}\, (q^{-1/2}z)^k,
\ee
are the polynomials considered by Pastro \cite{Pa} and for which
the
biorthogonality relation reads
\be
\frac{1}{2\pi i}\int_{|z| = 1} p_m(z, t_1, t_2)
 \overline{p_n(z, \overline{t_2},
\overline{t_1})}\,\Omega(z)\,\frac{dz}{z} = \frac{(t_1t_2q; q)_n}
{(q; q)_n}\, q^{-n} \delta_{m,n}.
\ee
Al-Salam and Ismail \cite {Al:Is} used (4.8) and the generating
function
\be
\Summ p_n(z; a, b)t^n =
 \frac{(atzq^{1/2}, bt; q)_\infty}{(tzq^{-1/2}, t;q)_\infty}
\ee
to establish a $q$-beta integral and found the rational functions
 biorthognal to
its integrand. The interested reader is refered to \cite{Al:Is} for
details.

\bigskip

\setcounter{equation}{0}

\setcounter{section}{5}
\setcounter{thm}{0}

{\bf  5. The $q^{-1}$-Hermite Ladder}. When $q > 1$ in (1.5) the
polynomials
 $\{H_n(x|q)\}$ become orthogonal on the imaginary axis. The result
of replacing
 $x$ by $ix$  and $q$ by $1/q$  put the orthogonality on the real
line and the new $q$ is now
in $(0,1)$, \cite {As}.  Denote $(-i)^n H_n(ix|1/q)$ by $h_n(x|q)$.
In this
new notation the recurrence relation (1.5) and the initial
conditions (1.7)
 become
\be
h_{n+1}(x|q)= 2x h_n(x|q) -q^{-n}(1-q^n)h_{n-1}(x|q), n > 0,
\ee
\be
h_0(x|q)= 1, \quad h_1(x|q) = 2x.
\ee
The polynomials $\{h_n(x|q)\}$ are called the $q^{-1}$-Hermite
polynomials,
 \cite {Is:Ma2}. The corresponding moment problem is indeterminate.
Let
 ${\cal V}_q$ be the set of probability measures which solve the
problem. For
 any $\mu\in {\cal V}_q$ we have

\be
\int_{-\infty}^\infty h_m(x|q) h_n(x|q) \, d\mu(x)
 = q^{-n(n+1)/2}(q; q)_n\delta_{m,n}.
\ee
The $h_n's$ have the generating function, \cite {Is:Ma2},
\be
\Summ h_n(x|q)\frac{t^n}{(q; q)_n}q^{n(n-1)/2}
= (-te^{\xi}, te^{-\xi}; q)_\infty, \quad x = \sinh \xi,\; t,\xi\in
\mbox{\bf C}.
\ee
By Proposition 1.1 it is clear that (5.3) and (5.4) imply
\be
\int_{-\infty}^\infty (-t_1e^{\xi}, t_1e^{-\xi}, -t_2e^{\xi},
t_2e^{-\xi};
 q)_\infty d\mu(x) = (-t_1t_2/q; q)_\infty,\quad
t_1,t_2\in\mbox{\bf C},
 \mu\in{\cal V}_q
\ee
since the power series (1.4) has radius of convergence
$\rho=\infty$.
Incidentally the function
\be
\chi_t(x)=(-te^{\xi},te^{-\xi};q)_\infty=
(-t(\sqrt{x^2+1}+x), t(\sqrt{x^2+1}-x); q)_\infty
\ee
belongs to $L^2(\mu)$ for any $\mu\in{\cal V}_q$ and any
$t\in\mbox{\bf C}$.
Therefore the complex measure $\nu_{\mu}(t_1,t_2)$ defined by
\be
d\nu_{\mu}(x; t_1, t_2) := \frac{\chi_{t_1}(x)\chi_{t_2}(x)}
{(-t_1t_2/q; q)_\infty} d\mu(x),\quad \mu\in{\cal V}_q, t_1,t_2\in
 \mbox{\bf C}, t_1t_2\neq -q^{1-k}, k\geq 0
\ee
has total mass $1$, and it is non-negative if $t_1=\overline{t_2}$.

Note that $(-te^{\xi},te^{-\xi};q)_k$ is a polynomial of degree $k$
in
$x=\sinh \xi$ for each fixed $t\neq 0$. Since
\[
(-te^{\xi}/q^k, te^{-\xi}/q^k;q)_k\chi_t(x)=\chi_{t/q^k}(x),
\]
we see that the non-negative polynomial
$|(-te^{\xi}/q^k,te^{-\xi}/q^k;q)_k|^2$
of degree $2k$ is $\nu_{\mu}(t,\overline{t})$-integrable. This
implies that
$\nu_{\mu}(t,\overline{t})$ has moments of any order, and by the
Cauchy-Schwarz
inequality every polynomial is $\nu_{\mu}(t_1,t_2)$-integrable for
all
$t_1,t_2\in\mbox{\bf C}, \mu\in{\cal V}_q$.

Introducing the orthonormal polynomials
\be
\tilde{h}_n(x|q)=\frac{h_n(x|q)}{\sqrt{(q;q)_n}}q^{n(n+1)/4}
\ee
the $q$-Mehler formula, cf. \cite{Is:Ma2}, reads
\be
\sum_{n=0}^\infty\tilde{h}_n(\sinh\xi|q)\tilde{h}_n(\sinh\eta|q)z
^n=
\frac{(-zqe^{\xi+\eta},-zqe^{-\xi-\eta},zqe^{\xi-\eta},zqe^{-\xi+
\eta};q)_\infty}
{(z^2q;q)_\infty}
\ee
valid for $\xi,\eta\in\mbox{\bf C}, |z|<1/\sqrt{q}.$

Applying the Darboux method to (5.9) Ismail and Masson
\cite{Is:Ma2} found
 the asymptotic
 behaviour of $\tilde{h}_n(\sinh\eta|q)$, and from their result it
follows
  that
 $(\tilde{h}_n(\sinh\eta|q)z^n)\in l^2$ for
$|z|<q^{-1/4},\eta\in\mbox{\bf C}$.
In this case the right-hand side of (5.9) belongs to $L^2(\mu)$ as
a function
of $x=\sinh\xi$ and the formula is its orthogonal expansion.
Putting
$t=qze^{\eta}, s=-qze^{-\eta}$, we have $z^2=-stq^{-2}$, so if
$|st|<q^{3/2}$
we have $|z|<q^{-1/4}$ and the right-hand side of (5.9) becomes
$\chi_t(\sinh\xi)\chi_s(\sinh\xi)/(-st/q;q)_\infty$, which belongs
to
 $L^2(\mu)$.
Using this observation we can give a simple proof of the following
formula
from \cite{Is:Ma2}.

\begin{pro} Let $\mu\in {\cal V}_q$ and let $t_i\in\mbox{\bf C},
i=1,\cdots,4$
satisfy $|t_1t_3|,|t_2t_4|<q^{3/2}$. (This holds in particular if
$|t_i|<
q^{3/4}, i=1,\cdots,4)$. Then $\prod_{i=1}^4\chi_{t_i}\in L^1(\mu)$
and
\be
\int \prod_{i=1}^4\chi_{t_i}\,d\mu=\frac{\prod_{1\leq j<k\leq
4}(-t_jt_k/q;q)_\infty}
{(t_1t_2t_3t_4q^{-3};q)_\infty}.
\ee
\end{pro}

\noindent{\bf Proof}. We write
\[
qz_1e^{\eta_1}=t_1,\; qz_2e^{\eta_2}=t_2,\; -qz_1e^{-\eta_1}=t_3,\;
 -qz_2e^{\eta_2}=t_4,
 \]
 noting that $z_1^2=-t_1t_3q^{-2}, z_2^2=-t_2t_4q^{-2}$, so the
equations have
 solutions $z_i,\eta_i,i=1,2$ if $t_i\neq 0$ for $i=1,\cdots,4$.
We next apply Parseval's formula to the two
 $L^2(\mu)$-functions $\chi_{t_1}\chi_{t_3},\;\chi_{t_2}\chi_{t_4}$
and get
\[
\int \prod_{i=1}^4\chi_{t_i}\,d\mu=(-t_1t_3/q,-t_2t_4/q;q)_\infty
\sum_{n=0}^\infty\tilde{h}_n(\sinh\eta_1|q)\tilde{h}_n(\sinh\eta_
2|q)
(z_1z_2)^n ,
\]
which by the $q$-Mehler formula gives the right-hand side of
(5.10).

 If $t_1=0$ and $t_2t_3t_4\neq 0$
 we apply Parseval's formula to $\chi_{t_3}$ and
$\chi_{t_2}\chi_{t_4}$, and if
 two of the parameters are zero the formula reduces to (5.5).
 \quad $\Box$

We shall now look at orthogonal polynomials for the measures
$\nu_{\mu}(t_1,t_2).$
When $q > 1$ the Al-Salam-Chihara polynomials are
 orthogonal on $(-\infty, \infty)$ and their moment problem may be
  indeterminate \cite {As:Is2}, \cite {Al:Ch}, \cite {Ch:Is}. If
one replaces
   $q$ by $1/q$ in the Al-Salam-Chihara polynomials, they can be
    renormalized to polynomials $\{v_n(x; q, a, b, c)\}$ 
satisfying

\be
(1-q^{n+1})v_{n+1}(x; q, a, b, c) = (a-xq^n)v_n(x; q, a, b, c) -
 (b - cq^{n-1})v_{n-1}(x; q, a, b, c),
\ee
where $0<q<1$, and $a,b,c$ are complex parameters. We  now consider
the
 special  case

\be
u_n(x; t_1, t_2) = v_n(-2x; q, -(t_1+t_2)/q, t_1t_2q^{-2}, -1),
\ee
where $t_1,t_2$ are complex parameters. The corresponding monic
polynomials
$\hat{u}_n(x)$ satisfy the recurrence relation determined from
(5.11)
\be
x\hat{u}_n(x)=\hat{u}_{n+1}(x)+\frac12(t_1+t_2)q^{-n-1}\hat{u}_n(
x)+
\frac14(t_1t_2q^{-2n-1}+q^{-n})(1-q^n)\hat{u}_{n-1}(x),
\ee
so by Favard's theorem, cf. \cite{Ch}, $\{u_n(x;t_1,t_2)\}$ are
orthogonal with
respect to a complex measure $\alpha(t_1,t_2)$ if and only if
$t_1t_2\neq-q^{n+1},n\geq 1$, and with respect to a probability
measure
$\alpha(t_1,t_2)$ if and only if $t_2=\overline{t_1}\in\mbox{\bf
C}\setminus
\mbox{\bf R}$ or $t_1,t_2\in\mbox{\bf R}$ and $t_1t_2\geq 0$. In
the affirmative
case
\be
\int_{-\infty}^\infty   u_m(x; t_1, t_2) u_n(x; t_1, t_2)\, d\a(x;
t_1, t_2)
 = \frac{q^{n(n-3)/2}}{(q; q)_n}(-t_1t_2q^{-n-1}; q)_n\delta_{m,n}.
\ee

It follows from (3.77) in \cite{As:Is2} that the moment problem is
indeterminate
if $t_1=\overline{t_2}$. If $t_1, t_2$ are real, different and
$t_1t_2\geq 0$,
we can assume $|t_1|<|t_2|$ without loss of generality, and in this
case the moment
problem is indeterminate if and only if $|t_1/t_2|>q$.

The generating function (3.70) in \cite {As:Is2} takes the form
$(x=\sinh\xi)$
\be
\Summ u_n(x; t_1, t_2) t^n =
\frac{(-te^\xi, te^{-\xi}; q)_\infty}{(-t_1t/q, -t_2t/q;
q)_\infty}, \quad
 \mbox{for}
  \; |t| < \; \min\{q/|t_1|, q/|t_2|\}.
\ee

By the method used in \cite{As:Is2} we can derive formulas for
$u_n$ in the
following way: Use the $q$-binomial theorem to write the right-hand
side of
 (5.15)
as a product of two power series in $t$ and equate coefficients of
$t^n$ to
get
\be
u_n(x;t_1,t_2)=\sum_{k=0}^n\frac{(qe^{\xi}/t_1; q)_k}{(q;
q)_k}(-t_1/q)^k
\frac{(-qe^{-\xi}/t_2; q)_{n-k}}{(q; q)_{n-k}}(-t_2/q)^{n-k}.
\ee
Application of the identity (I.11) in \cite{Ga:Ra}  gives the
explicit
 representation
 \be
u_n(x;t_1,t_2)=(-t_2/q)^n\frac{(-qe^{-\xi}/t_2; q)_n}{(q; q)_n}
{}_{2}\phi_1\left(\left. \ba{c} q^{-n}, qe^{\xi}/t_1 \\
-t_2e^{\xi}/q^n \ea \right|\,q,\;-t_1e^{\xi} \right),
\ee
which by (III.8) in \cite{Ga:Ra} can be transformed to

\be
u_n(x;t_1,t_2)=\frac{(-q^2/(t_1t_2); q)_n}{(q; q)_n}(-t_2/q)^n
{}_{3}\phi_1\left(\left. \ba{c} q^{-n}, qe^{\xi}/t_1,
-qe^{-\xi}/t_1 \\
-q^2/(t_1t_2) \ea \right|\,q,\;(t_1/t_2)q^n \right).
\ee

Writing the ${}_{3}\phi_1$ as a finite sum and applying the formula
\[
(a; q)_k=(q^{1-k}/a; q)_k(-a)^kq^{k(k-1)/2},
\]
we see that (5.18) can be transformed to
\bea
u_n(x;t_1,t_2)=
(-1/t_1)^n\frac{(-t_1t_2/q^{n+1}; q)_n}{(q; q)_n}q^{n(n+1)/2}
\sum_{k=0}^n
\frac{(q^{-n},-t_1e^{\xi}/q^k,t_1e^{-\xi}/q^k; q)_k}
{(q,-t_1t_2/q^{k+1}; q)_k}q^{nk}.
\eea

By symmetry of $t_1,t_2$ a similar formular holds for $t_1$ and
$t_2$
interchanged.

\begin{thm} For $\mu\in{\cal V}_q$ and $t_1,t_2\in\mbox{\bf C}$
such that
$t_1t_2\neq -q^{n+1},\; n\geq 1$ the Al-Salam-Chihara polynomials
$\{u_n(x;t_1,t_2)\}$ are orthogonal with respect to the complex
measure
 $\nu_{\mu}(t_1,t_2)$ given by (5.7).
 \end{thm}

\noindent{\bf Proof}. It follows by (5.13) that
$\hat{u}_n(x;0,0)=2^{-n}h_n(x|q)$
so the assertion is clear for $t_1=t_2=0$.

Assume now that $t_1\neq 0$. By the three term recurrence relation
it
 suffices to prove that
\be
\int u_n(x;t_1,t_2)\,d\nu_{\mu}(x;t_1,t_2)=0 \quad\mbox{for}\quad
n\geq 1.
\ee

By (5.19) and (5.7) we get
\bea
\int u_n(x;t_1,t_2)\,d\nu_{\mu}(x;t_1,t_2)
\nonumber
\eea
\bea
=(-1/t_1)^n\frac{(-t_1t_2/q^{n+1}; q)_n}{(q; q)_n}q^{n(n+1)/2}
\sum_{k=0}^n q^{nk}\frac{(q^{-n}; q)_k}{(q; q)_k}
\int\frac{\chi_{t_1/q^k}(x)\chi_{t_2}(x)}
{(-t_1t_2/q^{k+1}; q)_\infty}\,d\mu(x).
\nonumber
\eea

By (5.5) the integral is 1, and the sum is equal to 
${}_1\phi_0(q^{-n};-;q,q^n)$, which is equal to 0
 for $n\geq 1$
by (II.4) in \cite{Ga:Ra}.  \quad $\Box$

In particular, if $t_2=\overline{t_1}$ then
$\nu_{\mu}(t_1,\overline{t_1})$
is a positive measure and the Al-Salam-Chihara moment problem
corresponding to $\{u_n(x;t_1,\overline{t_1})\}$ is indeterminate.
The set
\[
\{\nu_{\mu}(t_1,\overline{t_1})\;|\;\mu\in{\cal V}_q \}
\]
is a compact convex subset of the full set ${\cal C}(t_1)$ of
solutions to
the $\{u_n(x;t_1,\overline{t_1})\}$-moment problem.

If $t_1=t_2\in\mbox{\bf R}$ then
\[
\{\nu_{\mu}(t_1,\overline{t_1})\;|\;\mu\in{\cal V}_q \}\neq {\cal
C}(t_1)
\]
since the measures on the left can have no mass at the zeros of
 $\chi_{t_1}(x)$.

If $t_1=t\in (q,1)$ and $t_2=0$ then the Al-Salam-Chihara moment
problem is
determinate and the set
$\{\nu_{\mu}(t,0)\;|\;\mu\in{\cal V}_q \}$ contains exactly one
positive
 measure namely the one coming from $\mu\in{\cal V}_q$ being the
 $N$-extremal solution corresponding to the choice
$a = q/t$ in (6.27) and (6.30)  of \cite {Is:Ma2}, i.e.
$\mu$ is the discrete measure with mass $m_n$ at $x_n$ for
$n\in\mbox{\bf Z}$,
where
\[
x_n=\frac12\left(\frac{t}{q^{n+1}}-\frac{q^{n+1}}{t}\right)
\]
and
\[
m_n=\frac{(q/t)^{4n}q^{n(2n-1)}(1+q^{2n+2}/t^2)}{(-q^2/t^2,-t^2/q
,q; q)_\infty}.
\]

The function $\chi_t(x)$ vanishes for $x=x_n$ when $n<0$ and we get
\[
\nu_{\mu}(t,0)=\sum_{n=0}^\infty c_n\varepsilon_{x_n},
\]
where
\[
c_n=\frac{q^{3n(n+1)/2}(1+q^{2n+2}/t^2)(-q^2/t^2; q)_n}{t^{2n}(q;
q)_n
(-q^2/t^2; q)_\infty }.
\]

We now go back to (5.5) and integrate $1/(-t_1e^{\xi},
t_1e^{-\xi};q)_k$
 against the integrand in (5.5). Here again the attachment method
works
  and we see that
\be
\varphi_n(\sinh \xi; t_1, t_2) := {}_{3}\phi_2\left(\left. \ba{c}
q^{-n}, -t_1t_2q^{n-2}, 0 \\
- t_1e^\xi,\;  t_1e^{-\xi} \ea \right|\,q,\;q \right).
\ee
satisfies the biorthogonality relation
\be
\int_{-\infty}^\infty \varphi_m(x; t_1, t_2)\varphi_n(x; t_2, t_1)
\chi_{t_1}(x)\chi_{t_2}(x)\, d\mu (x)
\ee
\bea
\qquad = \frac{1+t_1t_2q^{n-2}}{1+t_1t_2q^{2n-2}} (-t_1t_2q^{n-1};
q)_\infty(q; q)_n q^{n(n-3)/2}\, (t_1t_2)^n\,\delta_{m,n}.
\nonumber
\eea
The biorthogonal rational functions (5.21) are the special case
 $t_3 = t_4 = 0$
of the biorthogonal rational functions
\be
\varphi_n (\sinh\,\xi ;t_1,t_2,t_3,t_4) :=\,_4\phi_3
\left(\begin{array}{cc} q^{-n},-t_1t_2q^{n-2},-
t_1t_3/q,-t_1t_4/q\\
 -t_1e^\xi ,t_1e^{-\xi}  ,t_1t_2t_3t_4q^{-3}\end{array}\;\vrule\;\;
  q,q\right).
\ee
of Ismail and Masson \cite {Is:Ma2}.  We have not been able to
apply
 a generating function technique to (5.21) because we have not been
able to
  find a suitable generating function for the rational functions
(5.21).

We now return to the Al-Salam-Chihara polynomials $\{u_n(x; t_1,
t_2)\}$
in the positive definite case and
reconsider the generating function (5.15). The radius of
convergence of (1.4)
 is $\rho=q^{3/2}/\sqrt{t_1t_2}$, and we get by Proposition 1.1,
(5.14)
  and the $q$-binomial theorem
\be
\int_{-\infty}^\infty \chi_{t_3}(x)\chi_{t_4}(x)
d\a(x; t_1, t_2)
= \frac{(-t_1t_3/q, -t_1t_4/q, -t_2t_3/q, -t_2t_4/q, -t_3t_4/q;
q)_\infty}
{(t_1t_2t_3t_4q^{-3}; q)_\infty},
\ee
valid for $|t_3|,|t_4|<\rho$.

Applying this to the measures $\a(t_1,t_2)=\nu_{\mu}(t_1,t_2)$, we
get a new
proof of (5.10), now under slightly different assumptions on
$t_1,\cdots,t_4$.

The attachment procedure works in this case, and we prove the
biorthogonality
 relation of \cite {Is:Ma2} under the same assumptions as in
Proposition 5.1:

\bea
\int^\infty_{-\infty}\varphi_m (x;t_1,t_2,t_3,t_4)\;
\varphi_n(x;t_2,t_1,t_3,t_4)
 \prod_{i=1}^4\chi_{t_i}(x)\,d\mu(x)
\eea
\bea
=\frac{1+t_1t_2q^{n-2}}{1+t_1t_2q^{2n-2}}\,\frac{(t_1t_2t_3t_4q^{
-3})^n
(q, -q^2/t_3t_4; q)_n(-t_1t_2q^{n-1}; q)_\infty}
{(t_1t_2t_3t_4q^{-3}; q)_n}
\frac{\prod_{1\leq j<k\leq 4}(-t_jt_k/q;q)_\infty}
{(t_1t_2t_3t_4q^{-3};q)_\infty}\delta_{m,n}.\nonumber
\eea

Mathematics Institute, Copenhagen University, Universitetsparken 5,
DK-2100 Copenhagen \O, Denmark.

Department of Mathematics, University of South Florida, Tampa,
Florida 33620-5700.

\end{document}